\documentclass{amsart}
\usepackage{amssymb,amsmath, latexsym, mathbbol}
\theoremstyle{plain}
\newtheorem{theorem}{Theorem}

\newtheorem{lemma}{Lemma}

\theoremstyle{definition}

\theoremstyle{remark}

\numberwithin{equation}{section}
\newtheorem{remark}{Remark}
\numberwithin{equation}{section}
\newcommand{\e}{\epsilon}

\newcommand{\R}{\mathbb R}
\newcommand{\N}{\mathbb N}

\newcommand{\C}{\mathbb C}

\newcommand{\Rn}{\mathbb R^n}

\begin{document}
%
%
%
%
%
\title[Null--Control and measurable sets]{Null--Control and measurable sets}
\author{J. Apraiz}
\address[J. Apraiz]{Universidad del Pa{\'\i}s Vasco/Euskal Herriko
Unibertsitatea\\Departamento de Matem\'atica Aplicada\\Escuela Universitaria Polit\'ecnica de Donostia-San Sebasti\'an\\Plaza de Europa 1\\20018 Donostia-San Sebasti\'an, Spain.}
\email{jone.apraiz@ehu.es}
\author{L. Escauriaza}
\address[L. Escauriaza]{Universidad del Pa{\'\i}s Vasco/Euskal Herriko
Unibertsitatea\\Dpto. de Matem\'aticas\\Apto. 644, 48080 Bilbao, Spain.}
\email{luis.escauriaza@ehu.es}
\thanks{The authors are supported  by MEC grant, MTM2004-03029.}
\begin{abstract}
We prove the interior and boundary null--controllability of some parabolic evolutions with controls acting over measurable sets.
\end{abstract}
\maketitle

\section{Introduction}\label{S:1}
The control for evolution equations aims to drive the solution to a prescribed state starting
from a certain initial condition. One acts on the equation through a source term, a so-called distributed
control, or through a boundary condition. To achieve general results one
wishes for the control to only act in part of the domain or its boundary and to have as much
latitude as possible in the choice of the control region: location, size, shape.

Here, we focus on the heat equation in a smooth and bounded domain $\Omega$
 in $\Rn$ for a time interval $(0, T)$, $T > 0$ and for a distributed control $f$ we consider
\begin{equation}\label{E: system}
 \begin{cases}
\triangle u- \partial_tu = f(x,t) \chi_\omega(x),\  &\text{in}\ \Omega\times (0,T),\\
u=0,\ &\text{on}\ \partial\Omega\times [0,T],\\
u(0)= u_0,\ &\text{in}\ \Omega.
\end{cases}
\end{equation}
Here, $\omega\subset\Omega$ is an interior control region. The null controllability of this equation, i.e., the existence for
any $u_0$ in $L^2(\Omega)$ of a control $f$ in $L^2(\omega\times (0,T))$ with 
\begin{equation}\label{E:cost of control}
\|f\|_{L^2(\omega\times (0,T))}\le N\|u_0\|_{L^2(\Omega)},
\end{equation}
 such that $u(T)=0$, was proved in \cite{G. LebeauL. Robbiano} by means of local Carleman estimates for the elliptic operator $\triangle+\partial_y^2$ over $\Omega\times \R$. A second approach based on global Carleman estimates for the backward parabolic operator $\triangle+\partial_t$ \cite{FursikovOImanuvilov}, also led to the null controllability of the heat equation. The first approach has been used for the treatment of time-independent parabolic operators associated to self-adjoint elliptic operators, while the second allows to address time-dependent non-selfadjoint parabolic operators and semi-linear evolutions.
 
 The method introduced in  \cite{G. LebeauL. Robbiano} was further extended to study thermoelasticity \cite{G. LebeauE. Zuazua}, thermoelastic
plates \cite{BenabdallahNaso} and semigroups generated by fractional orders of elliptic operators \cite{Miller}. It has also been
used to prove null controllability in the case of non smooth coefficients \cite{BenabdallahDermenjianRousseau, RousseauRobbiano}.  The
method of \cite{G. LebeauL. Robbiano} has also be extended to treat some non-selfadjoint cases, e.g. non symmetric systems \cite{L«eautaud} and all 1-dimensional time-independent parabolic equations \cite{AlessandriniEscauriaza}.
 
 In the above works, the control region $\omega$ is always assumed to contain an open ball. Also, the cost of controllability (the smallest constant $N$ found for the inequality \eqref{E:cost of control}) depends on this fact. The reason for these is that the main technique used in the arguments, Carleman inequalities, requires to construct suitable Carleman weights: a role for functions which requires smoothness (at least $C^2$) and to have the extreme values in proper regions associated to the control region $\omega$, the larger body $\Omega$ and possibly the value of $T>0$. The construction of such functions seems to be not possible, when $\omega$  does not contain a ball.
 
Motivated by these facts J.P. Puel and E. Zuazua raised the question wether the null controllability of the heat equation is possible when the control region is a measurable set. A positive partial answer to this question was explained by the second author at the June 2008 meeting {\it Control of Physical Systems and Partial Differential Equations} held at the Institute Henri Poincar\'e. Here, we give a formal account of the results.
\begin{theorem}\label{T:1} Let $n\ge 2$. Then, $\triangle-\partial_t$ is null-controllable at all positive times, with distributed controls acting over a measurable set $\omega\subset\Omega$ with positive Lebesgue measure, when
\begin{equation*}
\triangle = \nabla\cdot\left(\mathbf A(x)\nabla\,\cdot\,\right)+V(x),
\end{equation*}
is a self-adjoint elliptic operator, the coefficients matrix $\mathbf A$ is smooth in $\overline\Omega$, $V$ is bounded in $\Omega$ and both are real-analytic in an open neighborhood of $\omega$. The same holds when $n=1$, 
\begin{equation*}
\triangle = \frac 1{\rho(x)}\left[\partial_x\left(a(x)\partial_x\,\right)+b(x)\partial_x\,+c(x)\right]
\end{equation*}
 and $a$, $b$, $c$ and $\rho$ are measurable functions in $\Omega=(0,1)$.
\end{theorem} 
 In regard to boundary null controllability, i.e., the existence for
any $u_0$ in $L^2(\Omega)$ of a control $h$ in $L^2(\gamma\times (0,T))$ with 
\begin{equation}\label{E:cost of control2}
\|h\|_{L^2(\gamma\times (0,T))}\le N\|u_0\|_{L^2(\Omega)},
\end{equation}
such that the solution to
 \begin{equation}\label{E: system2}
 \begin{cases}
\triangle u- \partial_tu = 0,\  &\text{in}\ \Omega\times (0,T),\\
u=h(x,t)\chi_\gamma(x),\ &\text{on}\ \partial\Omega\times [0,T],\\
u(0)= u_0,\ &\text{in}\ \Omega,
\end{cases}
\end{equation}
verifies $u(T)\equiv 0$, we have the following result.
 \begin{theorem}\label{T:2} Let $n\ge 2$. Then, $\triangle-\partial_t$ is null-controllable at all times $T>0 $ with boundary controls acting over a measurable set $\gamma\subset\partial\Omega$ with positive surface measure when
\begin{equation*}
\triangle = \nabla\cdot\left(\mathbf A(x)\nabla\,\cdot\,\right)+V(x)
\end{equation*}
is a self-adjoint elliptic operator, the coefficients matrix $\mathbf A$ is smooth in $\overline\Omega$, $V$ is bounded in $\Omega$ and both are real-analytic in an open neighborhood of $\gamma$ in $\overline\Omega$.
\end{theorem} 
The results in Theorems \ref{T:1} and \ref{T:2} follow from a straightforward application of the linear construction of the control function for the systems \eqref{E: system} and \eqref{E: system2} developed in \cite{G. LebeauL. Robbiano} and the following observability inequality or propagation of smallness estimate established in \cite{Vessella} (See also \cite{Nadirashvili2} and \cite{Nadirashvili}).
\begin{theorem}\label{T:3}
 Assume that $f:B_{2R}\subset\Rn\longrightarrow\R$ is a real-analytic function verifying
\begin{equation}\label{E: condicion fundament}
|\partial^\alpha f(x)|\le \frac{M |\alpha|!}{(\rho R)^{|\alpha|}}\ ,\  \text{when}\ \ x\in B_{2R}, \ \alpha\in \N^n,
\end{equation}
for some $M>0$ and $0<\rho\le 1$ and $E\subset B_{\frac R2}$ is a measurable set with positive measure. Then, there are positive constants $N=N(\rho, |E|/|B_R|)$ and $\theta=\theta(\rho, |E|/|B_R|)$ such that 
\begin{equation*}
\|f\|_{L^\infty(B_R)}\le N\left(\text{\rlap |{$\int_{E}$}}\,|f|\,dx\right)^{\theta}M^{1-\theta}.
\end{equation*}
\end{theorem}
The experts will realize that the word {\it smooth} describing the regularity (away from the measurable set) of $\partial\Omega$, $A$ and $V$ in Theorems \ref{T:1} and \ref{T:2} can be replaced by $\partial\Omega$ is $C^2$, $A$ is Lipschitz in $\overline\Omega$ and $V$ is bounded (See \cite{FursikovOImanuvilov,G. LebeauL. Robbiano,G. LebeauE. Zuazua,RousseauRobbiano2}). In fact, the $C^2$ regularity of $\partial\Omega$ can be relaxed to require that either $\partial\Omega$ is $C^1$ or there is $\alpha\in (0,1]$ such that
\begin{equation*}
\left(P-Q\right)\cdot \nu(P)\ge -|P-Q|^{1+\alpha},\ \text{for all}\ P, Q\in \partial\Omega,
\end{equation*}
where $\nu(P)$ is the exterior unit normal vector to $\partial\Omega$. The later holds when $\Omega$ is either a convex domain, a polyhedron in $\Rn$ or when $\partial\Omega$ can be locally written as the graphs of Lipschitz functions which are the sum of a convex and a $C^{1,\alpha}$ function over $\R^{n-1}$.

To simplify the exposition and to show the strength of Theorem \ref{T:3}, we give the proof of Theorems \ref{T:1} and \ref{T:2} under the assumptions that $\partial\Omega$, $A$ and $V$ are globally real analytic. We do it because it makes clear how the construction algorithm of the control function in \cite{G. LebeauL. Robbiano} and Theorem \ref{T:3} can also be applied to prove the {\it interior} null-controllability (Theorem \ref{T:1}) for other parabolic evolutions whose corresponding observability or spectral inequalities (suitable Carleman inequalities) are otherwise  unknown. Examples of these parabolic evolutions are the ones associated  to selfadjoint elliptic systems with unknowns $\mathbf u=(u^1,\dots,u^m)$, 
\begin{equation*}
L_\alpha \mathbf u=\partial_i(a^{\alpha\beta}_{ij}(x)\partial_ju^{\beta}),\ \alpha=1,\dots,m
\end{equation*}
with $a^{\alpha\beta}_{ij}=a^{\beta\alpha}_{ji}$, for $\alpha, \beta= 1,\dots, m$, $i, j=1,\dots,n$, and with coefficients matrices verifying for some $\delta >0$ the strong ellipticity condition,
\begin{equation*}
\sum_{i, j, \alpha,\beta}a^{\alpha\beta}_{ij}(x)\xi^\alpha_i\xi^{\beta}_j\ge\delta \sum_{i,\alpha}|\xi_i^\alpha|^2,\ \text{when}\ \boldsymbol\xi\in \R^{nm}, x\in\Rn,
\end{equation*}
or the more general Legendre-Hadamard condition 
\begin{equation}\label{E:Hadamard2}
\sum_{i, j, \alpha,\beta}a^{\alpha\beta}_{ij}(x)\xi_i\xi_j\eta^\alpha\eta^{\beta}\ge\delta | \xi|^2 |\boldsymbol\eta|^2,\ \text{when}\ \xi\in \R^{n}, \boldsymbol\eta\in\R^m, x\in\Rn.
\end{equation}
We recall that the Lam\'e system of elasticity
\begin{equation*}\label{E:segundosistema}
\nabla\cdot\left(\mu(x)\left(\nabla \mathbf u+\nabla \mathbf u^t\right)\right) +\nabla\left(\lambda(x)\nabla\cdot \mathbf u\right),
\end{equation*}
with $\mu\ge\delta$, $\mu+\lambda \ge 0$ in $\Rn$, $m=n$ and $a^{\alpha\beta}_{ij}=\mu(\delta_{\alpha\beta}\delta_{ij}+\delta_{i\beta}\delta_{j\alpha})+\lambda\delta_{j\beta}\delta_{\alpha i}$,  are examples of systems verifying \eqref{E:Hadamard2}. Here, $a^{\alpha\beta}_{ij}$, $\mu$ and $\lambda$ can either be constants or real analytic functions on $\overline\Omega$.

It also makes clear that under such hypothesis one may replace the Carleman inequalities used in the literature, to prove the observability or propagation of smallness inequalities necessary in the process of applying the construction algorithm in \cite{G. LebeauL. Robbiano,G. LebeauE. Zuazua}, by the simpler application of Theorem \ref{T:3}. Of course, it has the drawback that it requires more smoothness on the operators and the boundary of $\Omega$ but on the contrary, one can handle with Theorem \ref{T:3} and the construction methods in \cite{G. LebeauL. Robbiano,G. LebeauE. Zuazua} the null-controllability of other parabolic evolutions with internal controls: like the second order evolutions explained above or for higher order evolutions as \begin{equation*}
\partial_tu+(-1)^m\triangle^mu,\ m=2,\dots,
\end{equation*}
with Dirichlet boundary conditions  on $\partial\Omega$, $u=\nabla u=\dots=\nabla^{m-1}u=0$.

In section \ref{S:2} we give the proofs of Theorems \ref{T:1} and \ref{T:2}. We explain how to extend Theorem \ref{T:1} to the evolutions \eqref{E: system3} and \eqref{E:biharmonic} in Remark \ref{R:1}, while the problems we have to extend Theorem \ref{T:2} to these evolutions are explained in Remark \ref{R:2}, for the simpler case of parabolic systems. For the sake of completeness, we include a proof of Theorem \ref{T:3} in section \ref{S:3}. It is built with ideas taken from \cite{Malinnikova}, \cite{Nadirashvili2} and  \cite{Vessella}. 

\section{Proof of Theorems \ref{T:1} and \ref{T:2}}\label{S:2}
We begin by setting up the formal hypothesis: first we assume there is $0<\delta\le1$ such that
\begin{equation*}
\delta |\xi|^2\le\mathbf A(x)\xi\cdot\xi\le\delta^{-1} |\xi|^2,\ \text{for all}\  x\in\Omega,\ \xi\in\Rn,
\end{equation*}
$\Omega$ is a bounded open set in $\Rn$, $n\ge 2$, with a real analytic boundary and $\mathbf A$, $V$ are real analytic in $\overline\Omega$, i.e., there are $r>0$ and $0<\delta\le 1$ such that
\begin{equation*}
|\partial^{\alpha}\mathbf A(x)|+|\partial^\alpha V(x)|\le |\alpha|! \delta^{-|\alpha|-1}\, ,\ \text{when}\ x\in\overline\Omega,\ \alpha\in\N^n,
\end{equation*}
and for each $x\in\partial\Omega$, there are a new coordinate system (where $x=0$) and a real analytic function $\varphi : B_{r}'\subset \R^{n-1}\longrightarrow\R$ verifying
\begin{equation}\label{E:descripcionfrontera}
\begin{aligned}
\varphi(0'&)=0,\ |\partial^{\alpha}\varphi(x')|\le |\alpha|! \delta^{-|\alpha|-1}\, ,\ \text{when}\ x'\in B_{r}',\ \alpha\in \N^{n-1},\\
&B_{r}\cap\Omega=B_r\cap\{(x',x_n): x'\in B_{r}',\  x_n>\varphi(x')\},\\
&B_{r}\cap\partial\Omega=B_r\cap\{(x',x_n): x'\in B_{r}',\ x_n=\varphi(x')\}.
\end{aligned}
\end{equation}

When $E$ is a measurable set, $|E|$ will denote its Lebesgue or surface measure.

\begin{proof}[Proof of Theorem \ref{T:1}] We may assume that the eigenvalues with zero Dirichlet condition for $\triangle= \nabla\cdot\left(\mathbf A(x)\nabla\,\cdot\,\right)+V(x)$ on $\Omega$ are all positive, $0<\omega_1^2<\omega_2^2\le \omega_3^2\le \dots\le \omega_j^2\le \dots$ and $\{e_j\}$ denotes the sequence of $L^2(\Omega)$-normalized eigenfunctions,
\begin{equation*}
\begin{cases}
\triangle e_j+\omega_j^2e_j=0,\ &\text{in}\ \Omega,\\
e_j=0,\ &\text{in}\ \partial\Omega.
\end{cases}
\end{equation*}

When $\omega\subset\Omega$ is measurable with positive Lebesgue measure, the method in \cite{G. LebeauL. Robbiano} shows that one find and $L^2(\omega\times (0,T))$ control function $f$ verifying \eqref{E:cost of control} for the system \eqref{E: system}, provided there is $N=N(|\omega|,\Omega, r,\delta)$, such that the inequality
\begin{equation}\label{E: desigualdad}
\sum_{\omega_j\le\mu}a_j^2+b_j^2\le e^{N\mu}\int\int_{\omega\times [\frac 14,\frac 34]}|\sum_{\omega_j\le\mu}\left(a_je^{\omega_jy}+b_je^{-\omega_jy}\right)e_j|^2\,dxdy,
\end{equation}
holds for $\mu\ge\omega_1$ and all sequences $a_1, a_2,\dots$ and $b_1, b_2, \dots$ Let then,
\begin{equation}\label{E: unafuncionsimpatica}
u(x,y)=\sum_{\omega_j\le\mu}\left(a_je^{\omega_jy}+b_je^{-\omega_jy}\right)e_j,
\end{equation}
it satisfies, $\triangle u+\partial^2_yu=0$, in $\Omega\times\R$, $u=0$ on the lateral boundary of $\Omega\times\R$ and $u$ is real analytic in $\overline\Omega\times\R$. Moreover, given $(x_0,y_0)$ in $\overline\Omega\times\R$ and $R\le 1$, there are $N=N(r, \delta)$ and $\rho=\rho(r, \delta)$ such that
\begin{equation}\label{E: fundamental}
\|\partial^\alpha_x\partial_y^\beta u\|_{L^\infty(B_{R}(x_0,y_0)\cap\Omega\times\R)}\le \frac{N(|\alpha|+\beta)!}{(R\rho)^{|\alpha |+\beta}}\left(\text{\rlap |{$\int_{B_{2R}(x_0,y_0)\cap\Omega\times\R}$}}|u|^2\,dxdy\right)^{\frac 12},
\end{equation}
when $\alpha\in\N^n$ and $\beta\ge 1$. For the later see \cite[Chapter 5]{Morrey}, \cite[Chapter 3]{FJohn2}.

The orthonormality of $\{e_j\}$ in $\Omega$ and \eqref{E: fundamental} with $R=1$ imply
\begin{equation}\label{E: segundofundamental}
\|\partial_x^\alpha\partial_y^\beta u\|_{L^\infty(\Omega\times [-5,5])}\le e^{N\mu}(|\alpha|+\beta)!\rho^{-|\alpha |-\beta}\left(\sum_{\omega_j\le\mu}a_j^2+b_j^2\right)^{\frac 12},\ \text{for}\ \alpha\in\N^n, \beta\ge 0,
\end{equation}
and there is $C>0$ such that replacing the constants $N$ and $\rho$ in \eqref{E: segundofundamental} by $CN$ and $\rho/C$ respectively, $u$ has a real analytic extension to $\Omega_\rho\times [-4,4]$,
\begin{equation*}
\Omega_\rho=\{x\in\Rn: d(x,\Omega)\le\rho\},
\end{equation*}
with
\begin{equation}\label{E: segundofundamental2}
\|\partial_x^\alpha\partial_y^\beta u\|_{L^\infty(\Omega_\rho\times [-4,4])}\le M(|\alpha|+\beta)!(2\rho)^{-|\alpha |-\beta},\ \text{for}\ \alpha\in\N^n, \beta\ge 0,
\end{equation}
and
\begin{equation*}
M=e^{N\mu}\left(\sum_{\omega_j\le\mu}a_j^2+b_j^2\right)^{\frac 12}.
\end{equation*}
For $(x_0,y_0)$ in $\Omega\times [0,1]$ with $d(x_0,\partial\Omega)=\rho$, we have $B_{2\rho}(x_0,y_0)\subset\Omega_\rho\times [-4,4]$, and if we apply Theorem \ref{T:3} to the real analytic extension of $u$ in $B_{2\rho}(x_0,y_0)$ with $E=B_{\frac\rho 4}(x_0,y_0)$, \eqref{E: segundofundamental2} implies there is $0<\theta_1<1$ such that
\begin{equation*}
\|u\|_{L^\infty(B_\rho(x_0,y_0))}\le N\|u	\|_{L^\infty(B_{\frac\rho 4}(x_0,y_0))}^{\theta_1} M^{1-\theta_1},
\end{equation*}
From this and a suitable covering argument we get
\begin{equation}\label{E:primera}
\|u\|_{L^\infty(\Omega\times [0,1])}\le N\|u\|_{L^\infty(\Omega^\rho\times [-1,2])}^{\theta_1} M^{1-\theta_1},
\end{equation}
with
\begin{equation*}
\Omega^\rho=\{x\in\Omega : d(x,\partial\Omega)\ge \tfrac {3\rho}4\}.
\end{equation*} 
Thus, from Theorem 3 and without Carleman inequalities it is possible to bound except for the factor $M$ all the information related to $u$ over $\Omega\times [0,1]$ by the information on $u$ over $\Omega^\rho\times [-1,2]$, a region {\it located inside} $\Omega\times\R$.

We may also assume that $E=\omega\times (\frac {1-\rho}4,\frac {1+\rho}4)\subset B_{\frac\rho 2}(0,\frac 12)$ has positive measure inside $B_{\frac\rho 2}(0,\frac 12)$, $B_{8\rho}(0,\frac 12)\subset\Omega^\rho\times [0,1]$ and a second application of Theorem \ref{T:3} gives
\begin{equation}\label{E: segunda}
\|u\|_{L^\infty(B_\rho(0,\frac12))}\le N\|u\|_{L^2(\omega\times [\frac 14,\frac 34])}^{\theta_2} M^{1-\theta_2}.
\end{equation}
Proceeding as in \cite{G. LebeauL. Robbiano}, we use a covering of $\Omega^\rho\times [-1,2]$ and successive applications of Theorem \ref{T:3}, with $E$ being a centered-moving ball with fixed radius $R/2$ depending on $\rho$ and the geometry of $\Omega$, with the ball of the same center and radius $2R$ contained $\Omega_\rho\times [-4,4]$, and where in the last applications of Theorem \ref{T:3}, $E = B_\rho(0,\frac12)$. Thus, Theorem \ref{T:3} and \eqref{E: segundofundamental2} imply there is $0<\theta_3<1$ with
\begin{equation}\label{E:tercera}
\|u\|_{L^\infty(\Omega^\rho\times [-1,2])}\le N\|u\|_{L^\infty(B_\rho(0,\frac12))}^{\theta_3} M^{1-\theta_3}.
\end{equation}
The orthonormality of $\{e_j\}$ in $\Omega$ and the inequality
\begin{equation*}
e^{-\mu}\left(\frac{\sinh{\omega_1}}{\omega_1}-1\right)\left(a^2+b^2\right)\le\int_0^1\left(a e^{\omega y}+b e^{-\omega y}\right)^2dy,\ \text{when}\ \omega_1\le\omega\le\mu,\ a, b\in\R, 
\end{equation*}
give
\begin{equation}\label{E:cuarta}
\left(\sum_{\omega_j\le\mu}a_j^2+b_j^2\right)^{\frac 12}\le e^{N\mu}\|u\|_{L^2(\Omega\times [0,1])}\le e^{N\mu}|\Omega|^{\frac 12}\|u\|_{L^\infty(\Omega\times [0,1])} 
 \end{equation}
 and \eqref{E: desigualdad} follows from \eqref{E:cuarta}, \eqref{E:primera}, \eqref{E:tercera} and \eqref{E: segunda}.
 
In \cite{AlessandriniEscauriaza} it is shown  that the null-controllability of the system \eqref{E: system} over $\Omega= (0,1)$ with
\begin{equation*}
\triangle = \frac 1{\rho(x)}\left[\partial_x\left(a(x)\partial_x\,\right)+b(x)\partial_x\,+c(x)\right],
\end{equation*}
\begin{equation}\label{E: condici—n1}
\delta\le a(x),\ \rho (x)\le \delta^{-1},\ \ |b(x)|+|c(x)| \leq \delta^{-1}\ , \ \text{a.e. in} \ [0,1],
\end{equation}
is equivalent to the null-controllability of the system
\begin{equation}\label{E: parab—lica2}
\begin{cases}
\partial_{x}^2z-\rho(x)\partial_tz=f\chi_\omega, &\ 0<x<1\ ,\ 0<t<T,\\
z(0,t)=z(1,t)=0, &\ 0\le t\le T,\\
z(x,0)=z_0, & \ 0\le x\le 1,
\end{cases}
\end{equation}
where $\rho$ is a new function verifying \eqref{E: condici—n1} and $\omega$ a new measurable set in $[0,1]$ with positive measure. The later follows from the bilipschitz change of variables used in \cite{AlessandriniEscauriaza}. Let then, $0<\omega_1^2<\omega_2^2\le \omega_3^2\le \dots\le \omega_j^2\le \dots$ and $\{e_j\}$ denote the sequences of eigenvalues and $L^2(\Omega)$-normalized eigenfunctions verifying
\begin{equation*}
\begin{cases}
e_j''+\rho(x)\omega_j^2e_j=0, \ 0<x<1,\\
e_j(0)=e_j(1)=0.\\
\end{cases}
\end{equation*}
From \cite{G. LebeauL. Robbiano}, it suffices to show that \eqref{E: desigualdad} holds in order to find an interior null-control $f$ for \eqref{E: parab—lica2} verifying \eqref{E:cost of control}. Extend then $e_j$ and $\rho$ to
$[-1,1]$ by odd and even reflections respectively, and  to all $\R$ as periodic functions of period $2$. The extended $e_j$ is in $C^{1,1}(\R)$ and verifies $e_j''+\rho(x)\omega_je_j=0$ in $\R$, $j=1,2\dots$. As before, let $u$ be defined by \eqref{E: unafuncionsimpatica}, it verifies
\begin{equation*}
\partial^2_{x}u+\partial_y\left (\rho(x)\partial_{y}u\right)=0,\ \text{in}\ \R^2.
\end{equation*}
By Chebyshev's inequality, defining
\begin{equation*}
\left(\omega\times [\tfrac 14,\tfrac 34]\right)\setminus E=\{(x,y)\in \omega\times [\tfrac 14,\tfrac 34] : |u(x,y)|/2> \text{\rlap |{$\int_{\omega\times [\frac 14,\frac 34]}$}}|u|\,dxdy\},
\end{equation*}
we have
\begin{equation}\label{E:e–egir}
|E|\ge\tfrac 12|\omega\times [\tfrac 14,\tfrac 34]|\ \text{and}\ \|u\|_{L^\infty(E)}\le 2\  \text{\rlap |{$\int_{\omega\times [\frac 14,\frac 34]}$}}|u|\,dxdy.
\end{equation}
Set $u_\e(x,y)=u(\tfrac x\e,\tfrac y\e)$, it verifies
\begin{equation*}
\partial^2_{x}u_\e+\partial_y\left (\rho(x/\e)\partial_{y}u_\e\right)=0, \text{in}\ \R^2,
\end{equation*}
and let $v_\e$ be the ${\it stream}$ function of $u_\e$, i.e., the solution to
\begin{equation*}
\begin{cases}
\partial_xv_\e=-\rho(\tfrac x\e)\partial_yu_\e,\\
\partial_yv_\e=\partial_xu_\e,\\
v_\e(0)=0,
\end{cases}
\end{equation*}
Then, $f=u_\e+iv_\e$ is $(1/\delta)$-quasiregular, i.e.,
\begin{equation*}
f\in W^{1,2}_{\text{loc}}(\R^2),\ \partial_{\overline z}f=\nu(z)\partial_zf,\ |\nu(z)|\le \frac{1-\delta}{1+\delta},\ z\in\C,
\end{equation*}
and by the Ahlfors-Bers Representation Theorem \cite{AhlforsBers} (See also \cite{bn}
or \cite{bjs}), any $(1/\delta)$-quasiregular
mapping $f$ in $B_4$ can be written as
\begin{equation*}
f=F\circ\Psi,
\end{equation*}
where $F=U+iV$ is holomorphic in  $B_4$ and $\Psi : B_4\longrightarrow B_4$ is a $(1/\delta)$-quasiconformal mapping, i.e. a  $(1/\delta)$-quasiregular homeomorphism from $B_4$ onto $B_4$ verifying
\begin{equation}
\begin{aligned}\label{E: propiedadesconformal}
&\partial_{\overline z}\Psi=\nu(z)\partial_z\Psi,\ \Psi(0)=0, \Psi(4)=4,\\
N^{-1}|z_1-z_2|^{\frac 1{\alpha}}&\le  |\Psi(z_1)-\Psi(z_2)|\le
N|z_1-z_2|^{\alpha},\ \text{when}\ z_1, z_2\in B_4,
\end{aligned}
\end{equation}
for some $0<\alpha<1$ and $N\ge 1$ depending on $\delta$. Now, $\e E\subset B_{2\e}$, and from \eqref{E: propiedadesconformal}, $\Psi(\e E)\subset B_{C(2\e)^\alpha}$. Choose then $\epsilon$ so that $N(2\e)^\alpha=\frac 12$. Thus, $\Psi(\e E)\subset B_{\frac 12}$, $u_\e=U\circ \Psi$,
\begin{equation*}
\|U\|_{L^\infty(B_4)}=\|u\|_{L^\infty(B_{\frac 4\e})}
\end{equation*}
while the $L^\infty$ interior estimates for subsolutions of elliptic equations \cite[\S 8.6]{gt}, the periodicity and orthogonality of the eigenfunctions $e_j$ in $L^2([0,1],\rho\,dx)$, gives
\begin{equation*}
\|u\|_{L^\infty(B_{\frac 4\e})} \lesssim \|u\|_{L^2(B_{\frac 6\e})}\lesssim e^{6\mu /\e}\left(\sum_{\omega_j\le\mu}a_j^2+b_j^2\right)^{\frac 12}.
\end{equation*}
Thus, $U$ is harmonic in $B_4$,
\begin{equation*}
\|U\|_{L^\infty(B_4)}\le e^{N\mu}\left(\sum_{\omega_j\le\mu}a_j^2+b_j^2\right)^{\frac 12}
\end{equation*}
and from \eqref{E:e–egir}
\begin{equation}\label{E:e–egir3}
\|U\|_{L^\infty(\Psi(\e E))}\le 2\  \text{\rlap |{$\int_{\omega\times [\frac 14,\frac 34]}$}}|u|\,dxdy.
\end{equation}
All together, $U$ verifies the conditions in Theorem \ref{T:3} in $B_2$ with $R=1$, with the universal  constant $0<\rho\le 1$ associated to the quantitative analyticity over $B_2$ of  bounded harmonic functions in $B_4$. From \eqref{E:e–egir3}, \eqref{E: desigualdad} holds provided we can find a lower bound for the Lebesgue measure of $\Psi(\epsilon E)\subset B_{\frac 12}$. The lower bound follows from \eqref{E:e–egir} and the following rescaled version of \cite[Theorem 1]{Astala}: 

{\it Let $\Psi: B_4\longrightarrow B_4$ be a $(1/\delta)$-quasiconformal mapping with $\Psi(0)=0$ and $E\subset B_4$ be a measurable set. Then, there is $N=N(\delta)$ such that
\begin{equation*}
|E|^{\frac 1\delta}/N\le |\Psi(E)|\le N|E|^{\delta}.
\end{equation*}}
\noindent
\end{proof}
\begin{remark}\label{R:3}
Theorem  \ref{T:3} also implies the version of \eqref{E: desigualdad} appearing in \cite{G. LebeauE. Zuazua}. For if $\Omega$, $\mathbf A$ and $V$ are as above and
\begin{equation*}
u(x,y)=\sum_{\omega_j\le\mu}a_je^{\omega_jy}e_j(x),
\end{equation*}
$u$ verifies \eqref{E: fundamental} and
\begin{equation*}
\|\partial_x^\alpha u(\, .\,,0)\|_{L^\infty(\Omega)}\le M|\alpha|!(2\rho)^{-|\alpha |},\ \text{for}\ \alpha\in\N^n, \ \text{with}\ M=e^{N\mu}\left(\sum_{\omega_j\le\mu}a_j^2\right)^{\frac 12}
\end{equation*}
and $\rho$ as above. Thus, $u(\,.\, ,0)$ has an analytic extension to a $\rho$-neighborhood of $\Omega$, and after a finite number of applications of Theorem \ref{T:3} and a covering argument,
\begin{equation*}
\|u(\, .\, ,0)\|_{L^2(\Omega)}\le N\|u(\, .\, ,0)\|_{L^2(\omega)}^\theta M^{1-\theta}.
\end{equation*} 
In particular,
\begin{equation*}
\sum_{\omega_j\le\mu}a_j^2\le e^{N\mu}\int_{\omega}|\sum_{\omega_j\le\mu}a_je_j|^2\,dx,\ \text{with}\ N=N(|\omega|,\Omega,r,\delta).
\end{equation*}
\end{remark} 
\begin{proof}[Proof of Theorem \ref{T:2}] Let $u$ be defined by \eqref{E: unafuncionsimpatica}. From \cite{G. LebeauL. Robbiano}, one can find a boundary control $h$ verifying \eqref{E:cost of control2}, provided there is $N=N(|\gamma|,\Omega,r,\delta)$ such that the inequality

\begin{equation}\label{E: desigualdad5}
\sum_{\omega_j\le\mu}a_j^2+b_j^2\le e^{N\mu}\int\int_{\gamma\times [\frac 14,\frac 34]}|\sum_{\omega_j\le\mu}\left(a_je^{\omega_jy}+b_je^{-\omega_jy}\right)\frac{\,\partial e_j}{\partial n}|^2\,d\sigma dy,
\end{equation}
holds for $\mu\ge\omega_1$ and all sequences $a_1, a_2,\dots$ and $b_1, b_2, \dots$. Here, $\nu$, $\sigma$ and $\frac{\,\partial }{\partial n}$ denote respectively the exterior unit normal vector to $\Omega$, the surface measure on $\partial\Omega$ and the conormal derivative for $\partial^2_y+\triangle$ on $\partial\Omega\times\R$, $\frac{\,\partial e}{\partial n}=A\nabla_x e\cdot\nu$.  We may also assume that $0\in\partial\Omega$, $\gamma\subset B_{\frac \rho 2}\cap\partial\Omega$, where $B_{2\rho}\cap\partial\Omega$ is the region above the graph of a real analytic function $\varphi:B_{\rho}
'\subset\R^{n-1}\longrightarrow\R$, as in \eqref{E:descripcionfrontera}. From \cite[\S 3 (2)]{G. LebeauL. Robbiano}, there is $N$ such that
\begin{equation}\label{E: desigualdad8}
\sum_{\omega_j\le\mu}a_j^2+b_j^2\le e^{N\mu}\|\tfrac{\,\partial u}{\partial n}\|_{L^\infty(\partial\Omega\times [-1,2])},
\end{equation}
with
\begin{equation*}
\tfrac{\,\partial u}{\partial n}=\sum_{\omega_j\le\mu}\left(a_je^{\omega_jy}+b_je^{-\omega_jy}\right)\tfrac{\,\partial e_j}{\partial n}\, .
\end{equation*}
From \eqref{E: segundofundamental} and \eqref{E:descripcionfrontera}, there are $N=N(r,\delta)$ and $\rho=\rho(r,\delta)$ such that 
$h(x',y)=\frac{\,\partial u}{\partial n}(x',\varphi(x'),y)$ verifies
\begin{equation*}
\|\partial_{x'}^\alpha\partial_y^\beta h\|_{L^\infty(B_{2\rho}'\times [-4,4])}\le M\left(|\alpha|+\beta\right)!(2\rho)^{-|\alpha |-\beta},\ \text{for}\ \alpha\in\N^{n-1},\ \beta\in\N,
\end{equation*}
\begin{equation*}
M=e^{N\mu}\left(\sum_{\omega_j\le\mu}a_j^2+b_j^2\right)^{\frac 12},
\end{equation*}
when $B_{2\rho}\cap\partial\Omega$ is a coordinate chart of $\partial\Omega$ as in \eqref{E:descripcionfrontera}. This fact, a suitable covering argument of $\partial\Omega$ and the {\it three-spheres inequalities} associated to the {\it obvious} extension of Theorem \ref{T:3} for real analytic functions over a compact analytic surface in $\R^{n+1}$, imply there are $N=N(|\gamma|,r,\delta)$ and $\theta=\theta(|\gamma|,r,\delta)$ such that
\begin{equation}\label{E:yaterminando}
\|\tfrac{\,\partial u}{\partial n}\|_{L^\infty(\partial\Omega\times [-1,2])}\le N\|\tfrac{\,\partial u}{\partial n}\|_{L^2(\gamma\times [\frac 14,\frac 34])}^{\theta }M^{1-\theta}.
\end{equation}
Finally, \eqref{E: desigualdad5} follows from \eqref{E: desigualdad8}  and \eqref{E:yaterminando}.
\end{proof}
\begin{remark}\label{R:1} The extension of Theorem \ref{T:1} to the parabolic system
\begin{equation}\label{E: system3}
 \begin{cases}
\partial_i(a^{\alpha\beta}_{ij}\partial_je_k^\beta)-\partial_tu^\alpha=f^\alpha(x,t) \chi_\omega(x),\  &\text{in}\ \Omega_T,\ \alpha=1,\dots,m,\\
\mathbf u=0,\ &\text{on}\ \partial\Omega\times [0,T],\\
\mathbf u(0)= \mathbf u_0,\ &\text{in}\ \Omega.
\end{cases}
\end{equation}
with $\partial\Omega$ as in \eqref{E:descripcionfrontera}, $a^{\alpha\beta}_{ij}$ verifying \eqref{E:Hadamard2} and
\begin{equation}\label{E:quasiregularidad}
|\partial^{\gamma}a^{\alpha\beta}_{ij}(x)|\le  |\alpha|! \delta^{-|\alpha|-1}\, ,\ \text{when}\ x\in\overline\Omega,\ \gamma\in\N^n,
\end{equation}
for some $0<\delta\le 1$ is now obvious: the symmetry, coerciveness and compactness of the operator $L^2(\Omega)^m\longrightarrow W^{1,2}_0(\Omega)^m$, mapping $\mathbf f=(f^1,\dots,f^m)$ into the unique solution $\mathbf u=(u^1,\dots,u^m)$ to
\begin{equation*}
\begin{cases}
\partial_i(a^{\alpha\beta}_{ij}\partial_ju^\beta)-\Lambda u^\alpha=f^\alpha,\ &\text{in}\ \Omega,\ \alpha=1,\dots,m,\\
\mathbf u=0,\ &\text{in}\ \partial\Omega
\end{cases}
\end{equation*}
where $\Lambda>0$ is sufficiently large \cite[Prop. 2.1]{Giaquinta}, gives the existence of a complete system $\{\mathbf e_k\}$ in $L^2(\Omega)^m$, $\mathbf e_k=(e_k^1,\dots,e_k^m)$, of eigenfunctions verifying
\begin{equation*}
\begin{cases}
\partial_i(a^{\alpha\beta}_{ij}\partial_je_k^\beta)+\omega_k^2e_k^\alpha=0,\ &\text{in}\ \Omega,\ \alpha=1,\dots,m,\\
\mathbf e_k=0,\ &\text{in}\ \partial\Omega
\end{cases}
\end{equation*}
with eigenvalues $0\le\omega_1\le\dots\omega_k\le\dots$ and $\lim_{k\to +\infty}\omega_k=+\infty$. By separation of variables, the Green's matrix for the system \eqref{E: system3} over $\Omega\times\R$ is the $m\times m$ matrix
\begin{equation*}
\boldsymbol\Gamma(x,y,t-s)=\sum_{k=1}^{+\infty}e^{-\omega_k^2(t-s)}\mathbf e_k(x)\otimes \mathbf e_k(y).
\end{equation*}
Moreover, the interior and boundary regularity for the elliptic system $\partial^2_y+\partial_i(a^{\alpha\beta}_{ij}\partial_j\ )$ in $\Omega\times\R$, shows that \eqref{E: fundamental} holds for $\mathbf u$ as in \eqref{E: unafuncionsimpatica} but with $\mathbf e_k$ replacing $e_k$ (\cite{MorreyNirenberg}, \cite[Chapter II]{Giaquinta}). These and \cite{G. LebeauL. Robbiano} suffice to find a control function $\mathbf f$ for the system \eqref{E: system3} verifying \eqref{E:cost of control}. Furthermore, if you wish to get bounds on the regularity of $\mathbf f$, \cite{G. LebeauL. Robbiano} shows it suffices to know that, $\|e_k\|_{H^s(\Omega)}\le C_s(1+\omega_k)^s$, for $s\ge 0$, and that the number of $0\le\omega_k\le \mu$, is bounded by $N\mu^n$, when $\mu\ge 1$.The first holds from elliptic regularity and \eqref{E:quasiregularidad}, while the second follows from the Gaussian estimates verified by $\boldsymbol\Gamma$, i.e., there are $N$ and $\kappa$ \cite[Corollary 4.14]{ChoDongKim} such that
\begin{equation*}
|\boldsymbol \Gamma(x,y,t)|\le  N(1\wedge t)^{-\frac n2}e^{\Lambda t-\kappa |x-y|^2/t}, \text{for}\ x,y\in\Rn, t>0.
\end{equation*}
It implies,
\begin{equation*}
\int_{\Omega}\int_{\Omega}|\boldsymbol\Gamma(x,y,t)|^2\,dxdy=\sum_{k\ge 1}e^{-2\omega_k^2t}\le Ne^{2\Lambda t}|\Omega|t^{-\frac n2},
\end{equation*}
and suffices to choose $t=\frac 1{\mu^2}$. In particular, $\boldsymbol f\in C_0^\infty((0,T),C^\infty(\overline\Omega))$.

The null-controllability of the system
\begin{equation}\label{E:biharmonic}
\begin{cases}
\partial_tu+(-1)^m\triangle^mu=f(x,t)\chi\omega,\ &\text{in}\ \Omega\times (0,T],\\
u=\nabla u=\dots=\nabla^{m-1}u=0,\ &\text{in}\ \partial\Omega\times (0,T),\\
u(0)=u_0,\ &\text{in}\ \Omega,
\end{cases}
\end{equation}
$m\ge 2$, is better managed with the approach in \cite{G. LebeauE. Zuazua}. If $\{e_j\}$ and $0\le\omega_1^{2m}\le\dots\le\omega_k^{2m}\le\dots$ are the eigenvectors and eigenvalues for $\triangle^m$ in $W^{2,m}_0(\Omega)$,
\begin{equation*}
\begin{cases}
(-1)^m\triangle^me_j-\omega^{2m}_je_j=0,\ &\text{in}\ \Omega,\\
e_j=\nabla e_j=\dots=\nabla^{m-1}e_j=0,\ &\text{in}\ \partial\Omega,
\end{cases}
\end{equation*}
\begin{equation*}
u(x,y)=\sum_{w_j^m\le\mu} a_jX_j(y)e_j(x),
\end{equation*}
verifies, $\partial_y^{2m}u+\triangle^mu=0$ in $\Omega\times\R$, $u=\nabla u=\dots=\nabla^{m-1}u=0$ in $\partial\Omega\times \R$, when 
\begin{equation*}
X_j(y)=\begin{cases}
e^{\omega_jy},\ &\text{for}\ m\ \text{odd},\\
e^{\omega_je^{\tfrac{\pi i}{2m}}y},\  &\text{for}\ m\ \text{even}. 
\end{cases}
\end{equation*}
Again, $u$ verifies \eqref{E: fundamental} \cite{MorreyNirenberg} and from Theorem \ref{T:3} applied to $u(\, .\, ,0)$ as in Remark \ref{R:3},
\begin{equation*}
\sum_{\omega_j^m\le\mu}a_j^2\le e^{N\mu^{\frac 1m}}\int_{\omega}|\sum_{\omega_j^m\le\mu}a_je_j|^2\,dx,\ \text{with}\ N=N(|\omega|,\Omega,r,\delta, m).
\end{equation*}
From \cite{G. LebeauE. Zuazua}, the last inequality suffices to find a control $f$ in $L^2(0,T,L^2(\Omega))$ for \eqref{E:biharmonic}.
\end{remark}
\begin{remark}\label{R:2}
The proof of Theorem \ref{T:2} for the scalar case is based on the bound \eqref{E: desigualdad8}. In the literature it is obtained from \eqref{E:cuarta} and the interpolation inequality below, proved via Carleman inequalities: there are $N$ and $\theta$ such that
\begin{equation*}
\|u\|_{L^2(\Omega\times [0,1])}\le N\|\tfrac{\partial u}{\partial n}\|_{L^2(\partial\Omega\times [-1,2])}^\theta\|u\|_{L^2(\Omega\times [-3,3])}^{1-\theta},
\end{equation*}
holds whenever $u$ verifies $\triangle u+\partial^2_yu=0$ in $\Omega\times\R$ and $u=0$ on $\partial\Omega\times\R$. However, the authors are not aware wether  the corresponding interpolation inequality for solutions $\mathbf u$ of the elliptic system
\begin{equation*}
\begin{cases}
\partial_i(a^{\alpha\beta}_{ij}\partial_ju^\beta)+\partial^2_yu^\alpha=0,\ &\text{in}\ \Omega\times\R,\ \alpha=1,\dots,m,\\
\mathbf u=0,\ &\text{in}\ \partial\Omega\times\R,
\end{cases}
\end{equation*}
\begin{equation*}
\|\mathbf u\|_{L^2(\Omega\times [0,1])}\le N\|\tfrac{\partial \mathbf u}{\partial\mathbf n}\|_{L^2(\partial\Omega\times [-1,2])}^\theta\|\mathbf u\|_{L^2(\Omega\times [-3,3])}^{1-\theta},
\end{equation*}
holds. Here, $(\tfrac{\partial \mathbf u}{\partial\mathbf n})^\alpha=a^{\alpha\beta}_{ij}\partial_ju^\beta\nu_i$, $\alpha =1,\dots,m$, and for this reason we can not extend Theorem \ref{T:2} to parabolic systems. On the other hand, there is $\rho>0$ such that the mapping $\partial\Omega\times (0,\rho)\longrightarrow U_\rho$, $(Q,t)\rightsquigarrow Q+t\nu (Q)$, is an analytic diffeomorphism onto $U_\rho=\{x\in\Rn\setminus\overline\Omega : 0< d(x,\partial\Omega)<\rho\}$. Because the null-control of the parabolic system over $\overline\Omega\cup U_\rho$ and with controls acting over $\omega= U_{\frac{3\rho}4}\setminus U_{\frac \rho2}$ is possible, standard arguments show that the system
 \begin{equation*}
 \begin{cases}
\partial_i(a^{\alpha\beta}_{ij}\partial_je_k^\beta)-\partial_tu^\alpha=0,\  &\text{in}\ \Omega_T,\ \alpha=1,\dots,m,\\
\mathbf u=\mathbf g\ &\text{on}\ \partial\Omega\times [0,T],\\
\mathbf u(0)= \mathbf u_0,\ &\text{in}\ \Omega.
\end{cases}
\end{equation*}
can be null-controlled with controls $\mathbf g$ acting over the full lateral boundary of $\Omega$. 
\end{remark}
\section{Proof of Theorem \ref{T:3}}{\label{S:3}
First we recall Hadamard's three-circle theorem \cite{m} and prove two Lemmas before the proof of Theorem \ref{T:3}.
\begin{theorem}\label{T: 5}
 Let $F$ be a holomorphic function of a complex variable in the ball $B_{r_2}$. Then,
 the following is valid for $0<r_1\le r \le r_2$,
\begin{equation*}\label{E:Hadamard}
\|F\|_{L^\infty(B_r)}\le \|F\|_{L^\infty(B_{r_1})}^{\theta}\|F\|_{L^\infty(B_{r_2})}^{1-\theta},\ \theta= \frac{\log{\frac
{r_2}{r}}}{\log{\frac{r_2}{r_1}}}\,.
\end{equation*}
\end{theorem}
\begin{lemma}\label{L:1}
Let $f$ be holomorphic in $B_1$, $|f(z)|\le 1$ in $B_1$ and $E$ be a measurable set in $[-\frac 15,\frac 15]$. Then, there are $N=N(|E|)$ and $\gamma =\gamma(|E|)$ such that
\begin{equation*}
\|f\|_{L^\infty(B_{\frac 12})}\le N\|f\|_{L^\infty(E)}^{\gamma}
\end{equation*}
\end{lemma}
\begin{proof}
For $n\ge 1$, there are $n+1$ points with $-\frac 15\le x_0<x_1<\dots <x_n\le \frac 15$, with $x_i\in \overline E$, $i=0,\dots ,n$ and $x_i-x_{i-1}\ge\frac{|E|}{n+1}$, $i=1,\dots,n$. For example, $x_0=\inf{E}$, $x_i=\inf\left(E\cap [x_{i-1}+\frac{|E|}{n+1},\frac 15]\right)$. Let
\begin{equation*}
P_n(z)=\sum_{i=0}^n f(x_i)\,\frac{\prod_{j\neq i}(z-x_j)}{\,\,\prod_{j\neq i}(x_i-x_j)}
\end{equation*}
Then,
\begin{equation*}
|P_n(z)|\le \|f\|_{L^\infty(E)}|E|^{-n}\sum_{i=0}^n\frac{(n+1)^n}{i!(n-i)!}\le \|f\|_{L^\infty(E)}\left(\frac 3{|E|}\right)^n,\ \text{for}\ |z|\le \frac 12\, .
\end{equation*}
By Cauchy's formula,
\begin{equation*}
|f(z)-P_n(z)|=\left|\frac 1{2\pi i}\int_{|\xi|=1}\frac{f(\xi)(z-x_0)\dots(z-x_n)}{(\xi-z)(\xi-x_0)\dots(\xi-x_n)}\,d\xi\right|\le 2\left(\frac 78\right)^n,\ \text{for}\ |z|\le\frac 12.
\end{equation*}
The last two inequalities give
\begin{equation}\label{E: minimizar}
\|f\|_{L^\infty(B_{\frac 12})}\le \|f\|_{L^\infty(E)}\left(\frac 3{|E|}\right)^n+2\left(\frac 78\right)^n,\ \text{for all}\ n\ge 1,
\end{equation}
and the minimization in the $n$-variable of the right hand side of \eqref{E: minimizar} implies Lemma \ref{L:1}.
\end{proof}
\begin{lemma}\label{L:2}
Let $f$ be analytic in $[0,1]$, $E$ be a measurable set in $[0,1]$ and assume there are positive constants $M$ and $\rho$ such that
\begin{equation}\label{E: quantitativo}
|f^{(k)}(x)|\le M k!(2\rho)^{-k},\ \text{for}\ k\ge 0,\ x\in [0,1].
\end{equation}
Then, there are $N=N(\rho, |E|)$ and $\gamma=\gamma(\rho, |E|)$ such that
\begin{equation*}
\|f\|_{L^\infty([0,1])}\le N\|f\|_{L^\infty(E)}^{\gamma}M^{1-\gamma}.
\end{equation*}
\end{lemma}
\begin{proof} \eqref{E: quantitativo} implies that $f$ has a holomorphic extension to $D_{\rho}=\cup_{0\le x\le 1}B(x,\rho)$, with $|f|\le 2M$ in $D_{\rho}$. Write $[0,1]$ as a disjoint union of $\frac 5{2\rho}$ non-overlapping closed intervals of length $\frac {2\rho}5$. Among them there is at least one, $I=[x_0-\frac\delta 5,x_0+\frac \delta 5]$, such that $|E\cap I|\ge \frac{2\delta |E|}5$. Then, $g(z)=f(x_0+\delta z)/2M$ is holomorphic in $B_1$, $E_{x_0,\rho}=\rho^{-1}(E\cap I-x_0)$ is measurable in $[-\frac 15,\frac15]$ with measure bounded from below by $\frac{2|E|}5$, $\|g\|_{L^\infty(E_{x_0,\rho})}\le \|f\|_{L^\infty(E)}$ and applying Lemma \ref{L:1} to $g$
\begin{equation}\label{E:primerpaso}
\|f\|_{L^\infty(B_{\frac\rho 2}(x_0))}\le N\|f\|_{L^\infty(E)}^{\gamma}M^{1-\gamma},
\end{equation}
with $0\le x_0\le 1$. Finally, make successive applications of Hadamard's three-circle theorem (a finite number depending on $\rho$) with a suitable chain of three-circles of radius comparable to $\rho$ and with center at points $x$ in $[0,1]$ contained in $D_\rho$, while recalling that on the largest ball $|f|\le 2M$,  to get that 
\begin{equation}\label{E: sendafundament}
\|f\|_{L^\infty([0,1])}\le N\|f\|_{L^\infty(B_{\frac\rho 2}(x_0))}^\theta M^{1-\theta},\ \ \theta=\theta(\rho),
\end{equation}
and Lemma \ref{L:2} follows from \eqref{E: sendafundament} and \eqref{E:primerpaso}.
\end{proof}
\begin{proof}[Proof of Theorem \ref{T:3}] We may assume $R=1$. Let $x\in B_{\frac 12}$. Using spherical coordinates centered at $x$,
\begin{equation*}
|E|\le\int_{S^{n-1}}|\{t\in [0,1] : x+tz\in E\}|\,dz,
\end{equation*}
and there is at least one $z\in S^{n-1}$ with $|\{t\in [0,1] : x+tz\in E\}|\ge |E|/(2\omega_n)$, with $\omega_n$ the surface measure on $S^{n-1}$. Set $\varphi(t)=f(x+tz)$. From \eqref{E: condicion fundament}, $\varphi$ satisfies \eqref{E: quantitativo}, $\|\varphi\|_{L^\infty(E_z)}\le \|f\|_{L^\infty(E)}$ and Lemma \ref{L:2} gives
\begin{equation}\label{E: casicasi}
\|f\|_{L^\infty(B_{\frac 12})}\le N\|f\|_{L^\infty(E)}^{\gamma}M^{1-\gamma}.
\end{equation}
Finally, setting
\begin{equation*}
\widetilde E=\{x\in E: |f(x)|/2\le \text{\rlap |{$\int_{E}$}}\,|f|\,dx\},
\end{equation*}
Chebyshev's inequality shows that
\begin{equation*}
|\widetilde E|\ge |E|/2\ ,\ \|f\|_{L^\infty(\widetilde E)}\le 2\,\text{\rlap |{$\int_{E}$}}\,|f|\,dx,
\end{equation*}
and Theorem \ref{T:3} follows after replacing $E$ by $\widetilde E$ in \eqref{E: casicasi}.
\end{proof}
\medskip
\noindent {\it Acknowledgement}: The authors wish to thank S. Vessella for sharing his results in \cite{Vessella}.


\end{document}